\documentclass[a4paper]{article}
\usepackage{amsfonts}
\usepackage{amssymb}
\usepackage{amsmath}
\usepackage{amscd}
\usepackage{latexsym}
\usepackage{epsfig}
\newtheorem{theorem}{Theorem}[section]

\newcommand{\proof}{\medskip \noindent {\bf Proof. \ \ }}
\newcommand{\qed}{\null\hfill $\Box\;\;$ \medskip}
\begin{document}

\parbox{1mm}

\begin{center}
{\bf {\sc \Large On generalized Flett's mean value
theorem}}
\end{center}

\vskip 12pt

\begin{center}
{\bf Jana MOLN\'AROV\'A} \\ Institute of Mathematics, Faculty of
Science, Pavol Jozef \v Saf\'arik University in Ko\v sice,
Jesenn\'a 5, 040~01 Ko\v sice, Slovakia, \\ {\it E-mail address:}
jana.molnarova88@gmail.com
\end{center}

\vskip 24pt

\hspace{5mm}\parbox[t]{10cm}{\fontsize{9pt}{0.1in}\selectfont\noindent{\bf
Abstract.} We present a new proof of generalized Flett's mean
value theorem due to Pawlikowska (from 1999) using only the
original Flett's mean value theorem. Also, a Trahan-type condition
is established in general case.} \vskip 12pt

\hspace{5mm}\parbox[t]{10cm}{\fontsize{9pt}{0.1in}\selectfont\noindent{\bf
Key words and phrases.} Flett's mean value theorem, real function,
differentiability, Taylor polynomial} \vskip 24pt

%%%%%%%%%%%%%%%%%%%%%%%%%%%%%%%%%%%%%%%%%%%%%%%%%%%%%%%%%%%%%%%%%%%%%%%%%%%%%%%%%%%%%%%%%%%%%%%
\section{Introduction}
%%%%%%%%%%%%%%%%%%%%%%%%%%%%%%%%%%%%%%%%%%%%%%%%%%%%%%%%%%%%%%%%%%%%%%%%%%%%%%%%%%%%%%%%%%%%%%%

Mean value theorems play an essential role in analysis. The
simplest form of the mean value theorem due to Rolle is
well-known.

\begin{theorem}[Rolle's mean value theorem]\label{Rolle}
If $f: \langle a, b\rangle \to \mathbb{R}$ is continuous on
$\langle a,b\rangle$, differentiable on $(a,b)$ and $f(a) = f(b)$,
then there exists a number $\eta\in (a, b)$ such that
$f'(\eta)=0.$
\end{theorem}

A geometric interpretation of Theorem~\ref{Rolle} states that if
the curve $y = f(x)$ has a tangent at each point in $(a, b)$ and
$f(a) = f(b)$, then there exists a point $\eta\in (a, b)$ such
that the tangent at $(\eta, f(\eta))$ is parallel to the $x$-axis.
One may ask a natural question: \textsl{What if we remove the
boundary condition $f(a) = f(b)$}? The answer is well-known as the
Lagrange's mean value theorem. For the sake of brevity put
$$\phantom{}_a^b\mathcal{K}\left(f^{(n)},g^{(n)}\right) = \frac{f^{(n)}(b)-f^{(n)}(a)}{g^{(n)}(b)-g^{(n)}(a)}, \quad
n\in\mathbb{N}\cup\{0\},$$ for functions $f,g$ defined on $\langle
a,b\rangle$ (for which the expression has a sense). If
$g^{(n)}(b)-g^{(n)}(a)=b-a$, we simply write
$\phantom{}_a^b\mathcal{K}\left(f^{(n)}\right)$.

\begin{theorem}[Lagrange's mean value theorem]\label{Lagrange}
If $f: \langle a, b\rangle \to \mathbb{R}$ is continuous on
$\langle a,b\rangle$ and differentiable on $(a,b)$, then there
exists a number $\eta\in (a, b)$ such that
$f'(\eta)=\phantom{}_a^b\mathcal{K}(f)$.
\end{theorem}

Clearly, Theorem~\ref{Lagrange} reduces to Theorem~\ref{Rolle} if
$f(a)=f(b)$. Geometrically, Theorem~\ref{Lagrange} states that
given a line $\ell$ joining two points on the graph of a
differentiable function $f$, namely $(a, f(a))$ and $(b, f(b))$,
then there exists a point $\eta\in (a, b)$ such that the tangent
at $(\eta, f(\eta))$ is parallel to the given line $\ell$.

In connection with Theorem~\ref{Rolle} the following question may
arise: \textsl{Are there changes if in Theorem~\ref{Rolle} the
hypothesis $f(a) = f(b)$ refers to higher-order derivatives?}
T.~M. Flett, see~\cite{Flett}, first proved in 1958 the following
answer to this question for $n=1$ which gives a variant of
Lagrange's mean value theorem with Rolle-type condition.

\begin{theorem}[Flett's mean value theorem]\label{Flett}
If $f: \langle a, b\rangle \to \mathbb{R}$ is a differentiable
function on $\langle a,b\rangle$ and $f'(a) = f'(b)$, then there
exists a number $\eta\in (a, b)$ such that
\begin{equation}\label{2}
f'(\eta)=\phantom{}_a^\eta\mathcal{K}(f).
\end{equation}
\end{theorem}

\begin{figure}
\begin{center}
\includegraphics{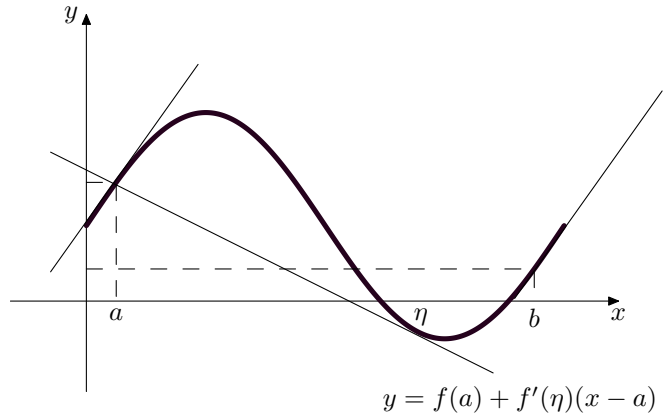}
\caption{Geometric interpretation of Flett's mean value theorem}
\label{obrFlett1}
\end{center}
\end{figure}

Flett's original proof, see~\cite{Flett}, uses
Theorem~\ref{Rolle}. A slightly different proof which uses
Fermat's theorem instead of Rolle's can be found in~\cite{RRA}.
There is a nice geometric interpretation of Theorem~\ref{Flett}:
if the curve $y = f(x)$ has a tangent at each point in $\langle a,
b\rangle$ and if the tangents at $(a, f(a))$ and $(b, f(b))$ are
parallel, then there exists a point $\eta \in (a,b)$ such that the
tangent at $(\eta, f(\eta))$ passes through the point $(a, f(a))$,
see Figure~\ref{obrFlett1}.

%Integral analogy of Theorem~\ref{Flett} was given by
%Wayment~\cite{Wayment}.
%\begin{theorem}[Wayment]\label{Wayment}
%If $f: \langle a, b\rangle \to \mathbb{R}$ is a continuous
%function on $\langle a,b\rangle$ and $f(a) = f(b)$, then there
%exists a number $\eta\in (a, b)$ such that
%$$f(\eta)(\eta - a) =  \int_a^\eta f(x)\,dx.$$
%\end{theorem}
%Wayment's proof, see also in~\cite{RS}, is unnecessarily
%complicated, because under the assumptions of
%Theorem~\ref{Wayment} it is enough to apply Flett's mean value
%theorem for function $F(x)=\int_a^x f(t)\,dt$.

Similarly as in the case of Rolle's theorem we may ask about
possibility to remove the boundary assumption $f'(a)=f'(b)$ in
Theorem~\ref{Flett}. As far as we know the first result of that
kind is given in book~\cite{RS}.

\begin{theorem}[Riedel-Sahoo]\label{thmRiedel-Sahoo}
If $f: \langle a, b\rangle \to \mathbb{R}$ is a differentiable
function on $\langle a,b\rangle$, then there exists a number
$\eta\in (a, b)$ such that
$$f'(\eta)=\phantom{}_a^\eta\mathcal{K}(f)+\phantom{}_a^b\mathcal{K}(f')\cdot\frac{\eta-a}{2}.$$
\end{theorem}

%Using this result applied for function $F(x)=\int_a^x f(t)\,dt$ it
%is possible to remove the condition $f(a) = f(b)$ in
%Theorem~\ref{Wayment} to obtain the following integral analogy of
%Theorem~\ref{thmRiedel-Sahoo}.
%\begin{theorem}
%If $f: \langle a, b\rangle \to \mathbb{R}$ is a continuous
%function on $\langle a,b\rangle$, then there exists a number
%$\eta\in (a, b)$ such that
%$$\int_a^\eta f(x)\,dx = f(\eta)(\eta - a) -\frac{1}{2}(\eta-a)^2 \frac{f(b)-f(a)}{b-a}.$$
%\end{theorem}

We point out that there are also other sufficient conditions
guaranteeing the existence of a point $\eta \in (a,b)$
satisfying~(\ref{2}). First such a condition was published in
Trahan's work~\cite{trahan}. An interesting idea is presented in
paper~\cite{tong} where the discrete and integral arithmetic mean
is used. We suppose that this idea may be further generalized for
the case of means studied e.g. in~\cite{hutnik, halhut, hutnik2}.

In recent years there has been renewed interest in Flett's mean
value theorem. Among the many other extensions and generalizations
of Theorem~\ref{Flett}, see e.g.~\cite{DPRS}, \cite{DRS},
\cite{JS}, \cite{PRS}, we focus on that of Iwona
Pawlikovska~\cite{Pawlikowska} solving the question of Zsolt Pales
raised at the 35-th International Symposium on Functional
Equations held in Graz in 1997.

\begin{theorem}[Pawlikowska]\label{thmPawlikowska}
Let $f$ be $n$-times differentiable on $\langle a,b\rangle$ and
$f^{(n)}(a)=f^{(n)}(b)$. Then there exists $\eta\in(a,b)$ such
that
\begin{equation}\label{genFlett}
f(\eta)-f(a)=\sum_{i=1}^n \frac{(-1)^{i+1}}{i!} (\eta-a)^i
f^{(i)}(\eta).
\end{equation}
\end{theorem}

Observe that the Pawlikowska's theorem has a close relationship
with the $n$-th Taylor polynomial of $f$. Indeed, for
$$T_n(f,x_0)(x) = f(x_0)+\frac{f'(x_0)}{1!}(x-x_0)+\dots +
\frac{f^{(n)}(x_0)}{n!}(x-x_0)^n$$ the Pawlikowska's theorem has
the following very easy form $f(a)=T_n(f,\eta)(a)$.

Pawlikowska's proof follows up the original idea of Flett,
see~\cite{Flett}, considering the auxiliary function
\begin{displaymath} G_f(x) = \left \{
\begin{array}{ll}
g^{(n-1)}(x), & x\in (a,b\rangle \\
\frac{1}{n}f^{(n)}(a), & x=a
\end{array} \right.
\end{displaymath} where $g(x)=\phantom{}_a^x\mathcal{K}(f)$ for $x\in
(a,b\rangle$ and using Theorem~\ref{Rolle}. In what follows we
provide a different proof of Theorem~\ref{thmPawlikowska} which
uses only iterations of an appropriate auxiliary function and
Theorem~\ref{Flett}. In Section~\ref{sectiontrahan} we give a
general version of Trahan condition, cf.~\cite{trahan} under which
Pawlikowska's theorem holds.

%%%%%%%%%%%%%%%%%%%%%%%%%%%%%%%%%%%%%%%%%%%%%%%%%%%%%%%%%%%%%%%%%%%
\section{New proof of Pawlikowska's theorem} \label{sectionnewproof}
%%%%%%%%%%%%%%%%%%%%%%%%%%%%%%%%%%%%%%%%%%%%%%%%%%%%%%%%%%%%%%%%%%%

The key tool in our proof consists in using the auxiliary function
$$\varphi_k(x)=x f^{(n-k+1)}(a)+\sum_{i=0}^{k} \frac{(-1)^{i+1}}{i!}(k-i)(x-a)^i
f^{(n-k+i)}(x), \quad k=1,2,\dots, n.$$ Running through all
indices $k=1, 2, \dots, n$ we show that its derivative fulfills
assumptions of Flett's mean value theorem and it implies the
validity of Flett's mean value theorem for $l$-th derivative of
$f$, where $l=n-1, n-2, \dots, 1$.

Indeed, for $k=1$ we have $\varphi_1(x)=-f^{(n-1)}(x)+x
f^{(n)}(a)$ and $\varphi_1'(x)=-f^{(n)}(x)+f^{(n)}(a)$. Clearly,
$\varphi_1'(a)=0=\varphi_1'(b)$, so applying the Flett's mean
value theorem for $\varphi_1$ on $\langle a,b\rangle$ there exists
$u_1\in(a,b)$ such that $\varphi_1'(u_1)(u_1-a) =
\varphi_1(u_1)-\varphi_1(a)$, i.e.
\begin{equation}\label{n-1}
f^{(n-1)}(u_1)-f^{(n-1)}(a) = (u_1-a)f^{(n)}(u_1).
\end{equation}
Then for $\varphi_2(x)=-2f^{(n-2)}(x)+(x-a) f^{(n-1)}(x)+x
f^{(n-1)}(a)$ we get
$$\varphi_2'(x)=-f^{(n-1)}(x)+(x-a)f^{(n)}(x)+f^{(n-1)}(a)$$ and $\varphi_2'(a)=0=\varphi_2'(u_1)$ by~(\ref{n-1}).
So, by Flett's mean value theorem for $\varphi_2$ on $\langle
a,u_1\rangle$ there exists $u_2\in (a,u_1)\subset (a,b)$ such that
$\varphi_2'(u_2)(u_2-a) = \varphi_2(u_2)-\varphi_2(a)$, which is
equivalent to $$f^{(n-2)}(u_2)-f^{(n-2)}(a) =
(u_2-a)f^{(n-1)}(u_2)-\frac{1}{2}(u_2-a)^2 f^{(n)}(u_2).$$
Continuing this way after $n-1$ steps, $n\geq 2$, there exists
$u_{n-1}\in (a,b)$ such that
\begin{equation}\label{1}
f'(u_{n-1})-f'(a) = \sum_{i=1}^{n-1} \frac{(-1)^{i+1}}{i!}
(u_{n-1}-a)^i f^{(i+1)}(u_{n-1}).
\end{equation} Considering the function $\varphi_n$ we get {\setlength\arraycolsep{2pt}
\begin{eqnarray*}
\varphi_n'(x) & = &
-f'(x)+f'(a)+\sum_{i=1}^{n-1}\frac{(-1)^{i+1}}{i!}(x-a)^i
f^{(i)}(x) \\ & = & f'(a)+\sum_{i=0}^{n-1}
\frac{(-1)^{i+1}}{i!}(x-a)^i f^{(i+1)}(x).
\end{eqnarray*}}Clearly, $\varphi_n'(a)=0=\varphi_n'(u_{n-1})$ by~(\ref{1}). Then by Flett's mean value theorem for
$\varphi_n$ on $\langle a,u_{n-1}\rangle$ there exists $\eta\in
(a,u_{n-1})\subset (a,b)$ such that
\begin{equation}\label{FMV1}
\varphi_n'(\eta)(\eta-a) = \varphi_n(\eta)-\varphi_n(a).
\end{equation} Since
$$\varphi_n'(\eta)(\eta-a) = (\eta-a)f'(a)+\sum_{i=1}^{n}
\frac{(-1)^{i}}{(i-1)!}(\eta-a)^i f^{(i)}(\eta)$$ and
$$\varphi_n(\eta)-\varphi_n(a) =
(\eta-a)f'(a)-n(f(\eta)-f(a))+\sum_{i=1}^{n}
\frac{(-1)^{i+1}}{i!}(n-i)(\eta-a)^i f^{(i)}(\eta),$$ the
equality~(\ref{FMV1}) yields {\setlength\arraycolsep{2pt}
\begin{eqnarray*}-n(f(\eta)-f(a)) & = & \sum_{i=1}^{n}
\frac{(-1)^{i}}{(i-1)!}(\eta-a)^i
f^{(i)}(\eta)\left(1+\frac{n-i}{i}\right) \\ & = & n
\sum_{i=1}^{n} \frac{(-1)^{i}}{i!}(\eta-a)^i f^{(i)}(\eta),
\end{eqnarray*}}which corresponds to~(\ref{genFlett}). \qed

It is also possible to state the result which no longer requires
any endpoint conditions. If we consider the auxiliary function
$$\psi_k(x)=\varphi_k(x)+\frac{(-1)^{k+1}}{(k+1)!}(x-a)^{k+1} \cdot
\phantom{}_a^b\mathcal{K}\left(f^{(n)}\right),$$ then the
analogous way as in the proof of Theorem~\ref{thmPawlikowska}
yields the following result also given in~\cite{Pawlikowska}
including Riedel-Sahoo's Theorem~\ref{thmRiedel-Sahoo} as a
special case ($n=1$).

\begin{theorem}
If $f: \langle a,b\rangle\to\mathbb{R}$ is $n$-times
differentiable on $\langle a,b\rangle$, then there exists
$\eta\in(a,b)$ such that
$$f(a) = T_n(f,\eta)(a)+\frac{(a-\eta)^{n+1}}{(n+1)!}\cdot
\phantom{}_a^b\mathcal{K}\left(f^{(n)}\right).$$
\end{theorem}

Note that the case $n=1$ is used to extend Flett's mean value
theorem for holomorphic functions, see~\cite{DPRS}. An easy
generalization of Pawlikowska's theorem involving two functions is
the following one.

\begin{theorem}\label{pawlikowskafg}
Let $f$, $g$ be $n$-times differentiable on $\langle a,b\rangle$
and $g^{(n)}(a)\neq g^{(n)}(b)$. Then there exists $\eta\in(a,b)$
such that
$$f(a)-T_n(f,\eta)(a)=
\phantom{}_a^b\mathcal{K}\left(f^{(n)},
g^{(n)}\right)\cdot[g(a)-T_n(g,\eta)(a)].$$
\end{theorem}

This may be verified applying the Pawlikowska's theorem to the
auxiliary function
$$h(x)=f(x)-\phantom{}_a^b\mathcal{K}\left(f^{(n)},
g^{(n)}\right)\cdot g(x), \quad x\in \langle a,b\rangle.$$ A
different proof will be presented in the following section.

%%%%%%%%%%%%%%%%%%%%%%%%%%%%%%%%%%%%%%%%%%%%%%%%%%%%%%%%%%%%%%%%%%%
\section{A Trahan-type condition}\label{sectiontrahan}
%%%%%%%%%%%%%%%%%%%%%%%%%%%%%%%%%%%%%%%%%%%%%%%%%%%%%%%%%%%%%%%%%%%

In~\cite{trahan} Trahan gave a sufficient condition for the
existence of a point $\eta \in (a,b)$ satisfying~(\ref{2}) under
the assumptions of differentiability of $f$ on $\langle
a,b\rangle$ and inequality
\begin{equation}\label{predpoklad}
(f'(b)-\phantom{}_a^b\mathcal{K}(f))\cdot
(f'(a)-\phantom{}_a^b\mathcal{K}(f))\geq 0.
\end{equation}

Modifying the Trahan's original proof using the Pawlikowska's
auxiliary function $G_f$ we are able to state the following
condition for validity~(\ref{genFlett}).

\begin{theorem}
Let $f$ be $n$-times differentiable on $\langle a,b\rangle$ and
$$\left(\frac{f^{(n)}(a)(a-b)^{n}}{n!}+M_f\right)
\left(\frac{f^{(n)}(b)(a-b)^{n}}{n!}+M_f\right)\geq 0,$$ where
$M_f=T_{n-1}(f,b)(a)-f(a)$. Then there exists $\eta\in(a,b)$
satisfying~(\ref{genFlett}).
\end{theorem}

\proof Since $G_f$ is continuous on $\langle a,b\rangle$ and
differentiable on $(a,b\rangle$ with $$G_f'(x) = g^{(n)}(x) =
\frac{(-1)^n
n!}{(x-a)^{n+1}}\left(f(x)-f(a)+\sum_{i=1}^{n}\frac{(-1)^i}{i!}(x-a)^i
f^{(i)}(x)\right),$$ for $x\in (a,b\rangle$,
see~\cite{Pawlikowska}, then {\setlength\arraycolsep{2pt}
\begin{eqnarray*} (G_f(b)-G_f(a))G_f'(b) & = &
\left(g^{(n-1)}(b)-\frac{1}{n}f^{(n)}(a)\right) g^{(n)}(b) \\ & =
& -\frac{n!(n-1)!}{(b-a)^{2n+1}}\left(\frac{f^{(n)}(a)(a-b)^{n}}{n!}+T_{n-1}(f,b)(a)-f(a)\right) \\
& \cdot &
\left(\frac{f^{(n)}(b)(a-b)^{n}}{n!}+T_{n-1}(f,b)(a)-f(a)\right)
\leq 0.
\end{eqnarray*}}According to Lemma~1 in~\cite{trahan} there exists $\eta\in (a,b)$
such that $G_f'(\eta)=0$ which corresponds to~(\ref{genFlett}).
\qed

Now we provide an alternative proof of Theorem~\ref{pawlikowskafg}
which does not use original Pawlikowska's theorem.

\vspace{0.3cm}\noindent\textbf{Proof of
Theorem~\ref{pawlikowskafg}.\,\,\,} For $x\in (a,b\rangle$ put
$\varphi(x)=\phantom{}_a^x\mathcal{K}(f)$ and
$\psi(x)=\phantom{}_a^x\mathcal{K}(g)$. Define the auxiliary
function $F$ as follows
\begin{equation}\label{fcia} F(x) =
\begin{cases}
\varphi^{(n-1)}(x)-\phantom{}_a^b\mathcal{K}\left(f^{(n)},
g^{(n)}\right)\cdot\psi^{(n-1)}(x),
& x\in(a,b\rangle \nonumber\\
\frac{1}{n}\left[f^{(n)}(a)-\phantom{}_a^b\mathcal{K}\left(f^{(n)},
g^{(n)}\right)\cdot g^{(n)}(a)\right], & x=a.\nonumber
\end{cases}
\end{equation}
Clearly, $F$ is continuous on $\langle a,b\rangle$, differentiable
on $(a,b\rangle$ and for $x\in(a,b\rangle$ there holds
{\setlength\arraycolsep{2pt}
\begin{eqnarray*}
F'(x) & = &
\varphi^{(n)}(x)-\phantom{}_a^b\mathcal{K}\left(f^{(n)},
g^{(n)}\right)\cdot\psi^{(n)}(x)
\\ & = & \frac{(-1)^n
n!}{(x-a)^{n+1}}\,\left[T_n(f,x)(a)-f(a)-\phantom{}_a^b\mathcal{K}\left(f^{(n)},
g^{(n)}\right)\cdot(T_n(g,x)(a)-g(a))\right].
\end{eqnarray*}}Then
$$F'(b)[F(b)-F(a)]=-\frac{n!(n-1)!}{(b-a)^{2n+1}}(F(b)-F(a))^{2}\leq 0,$$
and by Lemma~1 in~\cite{trahan} there exists $\eta\in(a,b)$ such
that $F'(\eta)=0$, i.e., $$f(a)-T_n(f,\eta)(a)=
\phantom{}_a^b\mathcal{K}\left(f^{(n)},
g^{(n)}\right)\cdot(g(a)-T_n(g,\eta)(a)).$$ \qed

\vskip 12pt

%%%%%%%%%%%%%%%%%%%%%%%%%%%%%%%%%%%%%%%%%%%%%%%%%%%%%%%%%%%%%%%%%%%%%%%%%%

\end{document}